\documentclass[11pt]{article}
\usepackage{mathrsfs}
\usepackage{amsfonts}
\usepackage{bbm}
\def\.{{\,\hskip-1pt\cdot\hskip-1pt\,}}

\def\b{{\mathbbm b}}

\def\qq{\quad{\rm and}\quad}

\input amssym.def
\input amssym.tex
\usepackage[all]{xy}
\usepackage{eufrak}
\usepackage{amscd}

\begin{document}
\begin{center}{\huge\bf
A remark on Rickard complexes}

\bigskip\bigskip{\bf\large Yuanyang Zhou}

\medskip{\scriptsize Department of Mathematics and
Statistics, Central China Normal University,
Wuhan, 430079, P.R. China

Email: zhouyuanyang@mail.ccnu.edu.cn}
\end{center}

\bigskip\noindent{\bf Abstract} In this paper, we characterize a Rickard complex, which induces a Rickard equivalence between the block algebras of a block $b$ and its Brauer correspondent and whose vertices have the same order as defect groups of the block $b$. The homology of such a Rickard complex vanishes at all degree but degree $q$, and the homology at degree $q$ induces a basic Morita equivalence between the block algebras in the sense of Puig.


\bigskip\bigskip\noindent{\bf 1.}\quad
In \cite{R}, J. Rickard exihibits a splendid Rickard complex, which induces a splendid Rickard equivalence between the block algebra of a $p$-block of a finite $p$-nilpotent group and the group algebra of its defect group, which is isomorphic to a Morita equivalence not induced by a $p$-permutation module. Then in \cite{H},  Harris and Linckelmann extend Rickard's technique and show a splendid Rickard  complex, which induces a splendid Rickard equivalence between the block algebras of a block for finite $p$-solvable groups with abelian defect groups and its Brauer correspondent, which is also isomophic to a Morita equivalence not induced by a $p$-permutation module. The two splendid Rickard complexes have vertex $ \Delta(Q)$ in terms of Puig (see Paragraph 7 below), where $Q$ is a defect group of the blocks and $\Delta(Q)$ is the diagonal subgroup of $Q\times Q$.  In this paper, we characterize a Rickard complex, which induces a Rickard equivalence between the block algebras of a block $b$ and its Brauer correspondent and whose vertices have the same order as defect groups of the block $b$.

\medskip\noindent{\bf 2.}\quad We recollect some notation in \cite{P3}.  Let $p$ be a prime number and let $k$ be an algebraically closed residue field $k$ of characteristic $p$.  From the point
of view of \cite{P3},  complexes are considered as
$\frak
D$-modules where, denoting by $\frak F$ the commutative
$k$-algebra of all the
$k$-valued functions on the set $\bf Z$ of all rational integers,
$\frak D$ is the
$k$-algebra containing $\frak{F}$ as a unitary $k$-subalgebra and an element $d$ such that
$$\frak{D}=\frak{F}\oplus \frak{F}d,\,\, d^2=0 \,\,\hbox{\rm and} \,\, df={\rm
sh}(f)d\neq 0 \,\,\hbox{\rm for any}\,\, f\in {\frak{F}-\{0\}}$$
where $\rm sh$ denotes the automorphism on the $k$-algebra $\frak{F}$
mapping $f\in \frak{F}$ onto the $k$-valued function sending $z\in \bf Z$ to
$f(z+1)\,;$ moreover, we denote by $s$ and $i_z\,,$ for any $z\in \bf Z\,,$ the
$k$-valued functions mapping $z'\in \bf Z$ on $(-1)^{z'}$ and
$\delta^{z'}_z$ respectively. Except for all the group algebras over
$\frak D\,,$ we
assume that all the modules and the algebras over $k$ are finitely generated. If $A$
is a $k$-algebra we denote by $A^*$ the group of invertible elements
of $A\,,$ and
by~$A^\circ$ the opposite $k$-algebra. Note that we have an
isomorphism ${\rm t}
\colon \frak D\cong \frak D^\circ$ mapping $f\in \frak{F}$ on the $k$-valued
function sending $z\in \bf Z$ to
$f(-z)\,,$ and $d$ on~$sd\,.$

\medskip\noindent{\bf 3.}\quad  A $\frak D$-interior algebra is a $k$-algebra
$A$ endowed with a
unitary $k$-algebra homomorphism $\varrho: \frak D\to A\,;$ for any $x,\, y\in \frak D$ and any $a\in A$, we write $x\cdot a\cdot y$ instead of $\varrho(x)a\varrho(y)$. Note that the isomorphism
${\rm t}\colon \frak D\cong \frak D^\circ$ then determines a $\frak
D$-interior
algebra structure for $A^\circ\,.$ Moreover, we have a $k $-algebra
homomorphism
$\frak D\to k$ mapping $f + f'd$ on $f(0)$ for any
$f,f'\in \frak F\,,$ so that any $k$-algebra admits a {\it
trivial\/} structure of
$\frak D$-interior algebra. The $\frak D$-interior algebra structure on $A$ induces a $\frak D$-module structure on $A$ by the equalities $$f(a)=\sum_{z,\, z'\in \bf Z}f(z)i_{z'}\cdot a\cdot i_{z'-z}\quad {\rm and} \quad d(a)=(d\cdot a-a\cdot d)\cdot s$$ for any $a\in A$ and any $f\in \frak F$. The $k$-algebra $A$ endowed with this $\frak D$-module is a $\frak D$-algebra in the sense of Puig (see \cite[11.2.4]{P3}).

\medskip\noindent
{\bf 4.} Let $G$ be a finite group; recall that a $k G$-interior algebra is a
$k$-algebra endowed with a unitary $k$-algebra homomorphism from
$k G\,.$ Similarly,
a $\frak DG$-interior algebra is a $k $-algebra $A$ endowed with a unitary
$k$-algebra homomorphism  $\rho\colon \frak DG\to A$ (but $A$ is always finitely
generated!); for any $x\in \frak DG$ and $a\in A\,,$ we write
$x\.a$ and $a\.x$
instead of $\rho(x)a$ and $a\rho(x)$ respectively.  If $A$ and $A'$ are ${\frak D}G$-interior
algebras, the
tensor product $A\otimes_k A'$ admits a ${\frak D}G$-interior algebra
structure given
by
$$f\.(a\otimes a') = \sum_{z,z'\in \bf Z} f(z +z')\,i_z\.a\otimes
i_{z'}\.a', \, \quad d\.(a\otimes a') = d\.a\otimes s\.a' + a\otimes d\.a'$$ and $g\.(a\otimes a')=g\.a\otimes g\.a'
$ for any $f\in \frak F$, any $g\in G$ and any $a,\, a'\in A$.
Here the first equality makes sense since in the sum above all but a finite number of
terms vanish and since we have ${\rm sh}(s) = -s\,.$
For any subgroup $H$ of
$G$, we denote by $A^H$ the centralizer of $\rho(H)$ in~$A\,;$
obviously $\rho(x)\in A^H$ for any $x\in \frak DC_G(H)$ and thus the
restriction of $\rho$ to $\frak DC_G(H)$ induces a
$\frak DC_G(H)$-interior algebra structure on $A^H$. Let $B$ and $C$ be two $k G$-interior algebras. A $k$-algebra homomorphism $f:B\rightarrow C$ is a $k G$-interior algebra homomorphism if $f$ preserves the $k G$-interior algebra structures on $B$ and $C$; furthermore, if $f$ is injective and $f(B)=f(1)Cf(1)$, then $f$ is a $k G$-interior algebra embedding.  Similarly ${\frak D}G$-interior algebra homomorphisms and ${\frak D}G$-interior algebra embeddings are defined.

\medskip\noindent
{\bf 5.} Let us denote by $\Bbb C_0 (A)$ the centralizer of the
image of $\frak D$
in $A\,;$  since the images of $\frak D$ and $G$ centralize each
other, $\Bbb C_0
(A)$ inherits a $k G$-interior algebra structure and, according to
the terminology
in \cite{P3}, the pointed groups, their inclusions, the local pointed
groups,\/ etc.
over the $\frak DG$-interior algebra~$A$ are nothing but the pointed
groups, their
inclusions, the local pointed groups, etc. over the $k G$-interior algebra
$\Bbb C_0 (A)\,.$ However, if $H_\beta$ is a pointed group over $A\,,$ so that
$\beta$ is a conjugacy class of primitive idempotents in $\Bbb C_0
(A)^H\,,$ the $k$-algebra $A_\beta = iAi$ for any
$i\in \beta$ inherits a $\frak DH$-interior
algebra structure mapping $y\in \frak DH$ on $y\.i = i\.y\,$; and the $k$-algebra $A_\beta$ is called an embedded algebra asociated with $H_\beta$.
For any
subgroup $H$
of $G\,,$ we call contractible any point contained in
the two-sided
ideal
$$\Bbb B_0 (A^H) = \Bbb C_0 (A)^H\cap \{d\.a + a\.d\mid a\in A^H\}$$
and we set $\Bbb H_0 (A^H) = \Bbb C_0 (A)^H/\Bbb B_0 (A^H)\,,$ which
still inherits
a $k C_G(H)$-interior algebra structure; whenever $\Bbb H_0 (A^G)
=\{0\}$ we say
that {\it $A$ is contractible\/}. It is clear that if $M$ is a $\frak
DG$-module
then ${\rm End}_k (M)$ is a ${\frak D}G$-interior algebra and we say
that $M$ is {\it
contractible\/} whenever ${\rm End}_k (M)$ is so \cite[Corollary 10.9]{P3};
moreover, we
say that $M$ is {\it $0$-split\/} if it is $\frak DG$-isomorphic to
the direct sum
of a contractible $\frak DG$-module and a $k G$-module endowed with
the trivial
$\frak D$-structure defined above.

\medskip\noindent
{\bf 6.} Let $G$ and $G'$ be two finite
groups and let $b$ and $b'$ be respective blocks of $G$ and $G'$.
Clearly the $k$-linear map $k G'\rightarrow k G'$ sending $x$ onto $x^{-1}$ for any $x\in G'$ is an opposite ring isomorphism. We denote by $b'^\circ$ the image of $b'$ through this opposite ring isomorphism. The $k$-linear map $k(G\times G')\rightarrow kG\otimes_k k G'$ sending $(x,\,y)$ onto $x\otimes y$ is a $k$-algebra isomorphism, through which we identify both sides so that $b\otimes b'^\circ$ is a block of $G\times G'$.
Let $\ddot M$ be an
indecomposable $\frak D(G\times G')$-module associated
with $b\otimes b'^\circ$ such that the restrictions of $\ddot M$ to $G\times \{1\}$ and
to $\{1\}\times G'$ are both projective. We denote by
$\ddot M^*$ the $k$-dual of~$\ddot M$ which, via the isomorphism
${\rm t}$ (see \cite[10.1.3]{P3}), still has a $\frak D (G\times G')$-module structure. Following
\cite[18.3.2]{P3}, we say that $\ddot M$ defines a Rickard equivalence between
$k Gb$
and~$k G'b'$ if, for
suitable contractible $\frak D(G\times G')$- and $\frak D(G'\times
G)$-modules
$C$~and~$C'\,,$ we have respective $\frak D(G\times G)$- and $\frak D(G'\times
G')$-module isomorphisms
$$\ddot M\otimes_{k G'} \ddot M^*\cong k Gb\oplus C\qq \ddot
M^*\otimes_{k G} \ddot M
\cong k G'b'\oplus C'$$
where $k Gb$ and $k G'b'$ have the trivial $\frak D$-interior
structure defined above.

\medskip\noindent
{\bf 7.} Let $\ddot P_{\ddot\gamma}$ be a maximal local pointed
group over the $\frak D(G\times
G')$-interior algebra
${\rm End}_k(\ddot M)$ or, equivalently, over the $k(G\times
G')$-interior algebra
$\Bbb C_0\big({\rm End}_k (\ddot M)\big)$. Then $\ddot P$ is a vertex of the $\frak D(G\times
G')$-module $\ddot M$ and a $\frak D\ddot P$-module $\ddot N=j\.\ddot M$ for
some $j\in
\ddot\gamma$ is a
source of the $\frak D(G\times
G')$-module $\ddot M$.
According to Theorem~18.8 in \cite{P3}, the images $P\subset G$
and $P'\subset G'$
of $\ddot P$ through the canonical projections
$\pi\colon G\times G'\to G$ and $\pi'\colon G\times G'\to G'$ are defect
groups of $b$ and $b'$ respectively.

\bigskip\noindent
{\bf Theorem 8.} {\it Let $G$ be a finite group, let $b$ be a block of
$G$ with defect group $P$ and let $b'$ be the Brauer correspondent
of $b$ in the normalizer $G'$ of $P$ in $G$. Let $\ddot M$ be a noncontractible
indecomosable ${\frak D}(G\times G')$-module inducing a Rickard equivalence between $k G b$ and $k G'b'$.  Let ${\frak D}\ddot P$-module $\ddot N$ be a source of the ${\frak D}(G\times G')$-module $\ddot M$.
Then the following are equivalent:

\medskip\noindent
{\bf 8.1.}
The groups $\ddot P$ and $P$ have the same order.

\medskip\noindent
{\bf 8.2.}
The homology of the ${\frak D}(G\times G')$-module $\ddot M$ vanishes at all degree but some degree q, the homology ${\Bbb H}_q(\ddot M)$ at degree $q$,
a $k(G\times G')$-module, induces a basic Morita equivalence
between $k G b$
and $k G' b'$ in the sense of Puig in \cite{P3}, and the ${\frak D}\ddot P$-module ${\rm End}_k(\ddot N)$ determined by the $\ddot P$-conjugation and the $\frak D$-algebra structure on ${\rm End}_k(\ddot N)$ is $0$-split.

\medskip\noindent In this case, the homology of the ${\frak D} \ddot P$-module $\ddot N$ vanishes at all degree but degree $q$, and the  homology ${\Bbb H}_q(\ddot N)$ at degree $q$, a $k \ddot P$-module, is a source  of the $k(G\times G')$-module ${\Bbb H}_q(\ddot M)$.}

\bigskip\noindent{\bf 9.}\quad We recall Brauer quotients and Brauer homomorphisms in \cite{T} and then prepare several lemmas.
Let $G$ be a finite group and $V$ be a $kG$-module. For
any subgroup $P$
of~$G$, we denote by $V^{P}$ the $k$-submodule of all $P$-fixed
elements of $V$,
by~$V(P)$ the Brauer quotient
$$ V(P)=V^{P}\Big/\sum_RV^{P}_{R}
\quad, $$
where $R$ runs over the set of all proper subgroups of $P$ and
$V^{P}_{R}$ is the image of the usual {relative trace map} ${\rm
Tr}^{P}_{R}: V^{R}\rightarrow V^{P}$, and by ${\rm Br}_P^V$ the
canonical surjective homomorphism $V^P \rightarrow V(P)$, which is
the so-called {\it Brauer homomorphism\/} associated to $P$ and $V$. Obviously,
the $kG$-module structure on~$V$ induces  $kN_G(P)$-module
structures on both $V^P$ and $V(P)\,,$ and ${\rm Br}_P^V$ is a
homomorphism of $k
N_G(P)$-modules. Let $A$ be a $k G$-interior algebra. We apply the Brauer quotient $V(P)$ to the $k G$-module $A$ induced by the $G$-conjugation, and then get the Brauer quotient $A(P)$. It is easily checked that both $A^P$ and $A(P)$ are $k$-algebras and that the Brauer homomorphism ${\rm Br}_P^A$ is a $k$-algebra homomorphism, whose kernel is the sum of all ideals $A^P_R$, where $R$ runs over the set of all proper subgroups of $P$.

\bigskip Let $Q$ be a $p$-group.  A $k Q$-interior algebra $A$ is a primitive $k Q$-interior algebra if $A^Q$ is a local algebra; furthermore, if $A(Q)\neq 0$, then $A$ is a local primitive $k Q$-interior algebra.

\bigskip \noindent {\bf Lemma 10.}\quad {\it Let $H$ be a finite group with a normal $p$-subgroup $Q$ and $T$ be a local primitive $kQ$-interior algebra. Then the kernel ${\rm Ker}({\rm Br}_Q^{T\otimes_k k H})$ is contained in the radical of $(T\otimes_k k H)^Q$.}

\medskip\noindent{\it Proof.}\quad Set $\bar H=H/Q$. The canonical homomorphism $H\rightarrow {\bar H}$ induces a $k$-algebra homomorphism $\varpi:k H\rightarrow k \bar H$. By tensorin both sides of the homomorphism $\varpi$ with $T$, we get a new $k$-algebra homomorphism $$1\otimes \varpi: T\otimes_k k H\rightarrow T\otimes_k  k\bar H$$ mapping $t\otimes a$ onto $t\otimes \varpi(a)$ for any $t\in T$ and $a\in k H$. Clearly $1\otimes \varpi$ is a $kQ$-interior algebra homomorphism and it induces a $k$-algebra homomorphism $$(T\otimes_k k H)^Q\rightarrow (T\otimes_k k\bar H)^Q=T^Q\otimes_k k\bar H,\leqno 10.1$$
which maps $(T\otimes_k k H)^Q_R$ into $T^Q_R\otimes_k k\bar H$ for any proper subgroup $R$ of $Q$. Since $T$ is a local primitive $kQ$-interior algebra, $T^Q_R$ is contianed in $J(T^Q)$ and thus the image of $(T\otimes_k k H)^Q_R$ through Homomorphism 10.1 is contained in the radical of the image of Homomorphism 10.1. But obviously the kernel of Homomorphism 10.1 is contained in the radical of $(T\otimes_k k H)^Q$. Therefore $(T\otimes_k k H)^Q_R$ is contained in the radical of $(T\otimes_k k H)^Q$. The proof is done.

\bigskip \noindent {\bf Lemma 11.}\quad {\it Let $H$ be a finite group with a normal $p$-subgroup $Q$ and $T$ be a local primitive $Q$-interior algebra. If $i$ is a primitive idempotent in $(k H)^Q$ such that ${\rm Br}_Q^{k H}(i)\neq 0$, then $1\otimes i$ is also a primitive idempotent in $(T\otimes_k k H)^Q$ such that ${\rm Br}_Q^{T\otimes_k k H}(1\otimes i)\neq 0$.}

\smallskip\noindent{\it Proof.}\quad Since ${\rm Br}_Q^T(1)$ is primitive in $T(Q)$ and ${\rm Br}_Q^{k H}(i)$ is primitive in $(k H)(Q)$,  ${\rm Br}_Q^T(1)\otimes {\rm Br}_Q^{k H}(i)$ is primitive in $T(Q)\otimes_k (k H)(Q)$.
 On the other hand, by \cite[Proposition 5.6]{P4}, there is a $k$-algebra isomorphism $$T(Q)\otimes_k (k H)(Q)\cong (T\otimes _k k H)(Q)$$ mapping ${\rm Br}_Q^T(t)\otimes{\rm Br}_Q^{k H}(a)$ onto ${\rm Br}_Q^{T\otimes_k k H}(t\otimes a)$.
In particular, this isomorphism maps ${\rm Br}_Q^T(1)\otimes{\rm Br}_Q^{k H}(i)$ onto ${\rm Br}_Q^{T\otimes_k k H}(1\otimes i)$.
Therefore ${\rm Br}_Q^{T\otimes_k k H}(1\otimes i)$
is primitive in $(T\otimes _k k H)(Q)$.
Since it follows from Lemma 10 that the kernel
${\rm Ker}({\rm Br}_Q^{T\otimes_k k H})$ is
contained in the radical of $(T\otimes_k k H)^Q$,
$1\otimes i$ is a primitive idempotent in
$(T\otimes_k k H)^Q$.

\bigskip \noindent {\bf Lemma 12.}\quad {\it Let $M$ be a
${\frak D}P$-module. Assume that the $k P$-interior algebra ${\Bbb H}_0({\rm End}_k(M))$ is a primitive $k P$-interior algebra. Then
the homology of $M$ vanishes at all degree but some degree
$q$ and there is a
$k P$-interior algebra isomorphism
${\Bbb H}_0({\rm End}_k(M))\cong {\rm End}_k({\Bbb H}_q(M))$.}

\medskip\noindent{\it Proof.}\quad By \cite[Theorem 3.1]{HS}, there is a short exact sequence of group homomorphisms $$0\rightarrow\Pi_q {\rm Ext}_k^1({\Bbb H}_q(M),\,{\Bbb H}_{q+1}(M))\rightarrow {\Bbb H}_0({\rm End}_k(M))\buildrel \xi \over{\rightarrow} \Pi_q{\rm End}_k({\Bbb H}_q(M))\rightarrow 0,$$ where $\xi$ is induced by the map ${\Bbb C}_0({\rm End}_k(M))\rightarrow \Pi_q{\rm End}_k({\Bbb H}_q(M))$ sending a chain map $f$
onto the induced family $(f_q)\in \Pi_q{\rm End}_k({\Bbb H}_q(M))$.
Since ${\rm Ext}_k^1({\Bbb H}_q(M),\,{\Bbb H}_{q+1}(M))=0\,\, $ for each $q$, the homomorphism $\xi$ is a group isomorphism.  Clearly the isomorphism $\xi$ preserves the composition of maps and the $k P$-interior algebra structures on ${\Bbb H}_0({\rm End}_k(M))$ and $\Pi_q{\rm End}_k({\Bbb H}_q(M))$; that is to say, $\xi$ is a $k P$-interior algebra isomorphism.
Since the $k P$-interior algebra ${\Bbb H}_0({\rm End}_k(M))$ is a primitive $k P$-interior algebra, the isomorphism $\xi$ forces that
the homology of $M$ vanishes at all degree but some degree
$q$; in particular, we have a $k P$-algebra isomorphism ${\Bbb H}_0({\rm End}_k(M))\cong {\rm End}_k({\Bbb H}_q(M))$. The proof is done.

\bigskip \noindent {\bf Remark.}\quad The differential $d$ on the $\frak D$-algebra ${\rm End}_k(M)$ (see Paragraph 3) is different from a differential on ${\rm End}_k(M)$ defined in \cite[Charpter V, 1.6]{HS}, but by replacing $d$ by $sd$, this difference will disappear.
Since $s$ does not affect the homology of ${\rm End}_k(M)$, \cite[Theorem 3.1]{HS} can be applied to the $\frak D$-algebra ${\rm End}_k(M)$.

\medskip\noindent{\bf 13.}\quad We begin to prove Theorm 8. We keep the notation in Theorem 8 and assume that Statement 8.1 holds.
By \cite[Theorem 18.8]{P3}, the images  $R$ and $R'$
of $\ddot P$ through the canonical projections
$\pi\colon G\times G'\to G$ and $\pi'\colon G\times G'\to G'$ are defect
groups of $b$ and $b'$ respectively. Since the orders of $P$ and $\ddot P$ are the same, the two projections $\pi$ and $\pi'$ induces group isomorphisms $\ddot P\cong R$ and $\ddot P\cong R'$, through which we identify $R$, $R'$ and $\ddot P$. Since $P$ is normal in $G'$, $P$ is the unique defect group of $b'$ and thus $R'$ is equal to $P$. In particular, we have $R=\ddot R=R'=P$.
By \cite[Theorem 18.8]{P3},
there are maximal local
pointed groups $P_\gamma$ on $k G b$, $P_{\gamma'}$ on
$k G' b' $ and $P_{\ddot\gamma}$ on ${\rm End}_k (\ddot N)\otimes_ k (k G')_
{\gamma'}$, such that we have a $k P$-interior algebra
isomorphism $$(k G)_\gamma\cong {\Bbb H}_0(({\rm End}_k (\ddot N)\otimes_ k (k G')_{\gamma'})_{\ddot\gamma})\leqno 13.1$$ and such that
the ${\frak D}({P}\times {P})$-module $({\rm End}_k (\ddot N)\otimes_ k (k G')_{\gamma'})_{\ddot\gamma}$ determined by the $\frak D$-algebra structure on ${\rm End}_k (\ddot N)$ and by the left and right multiplications of $P$ on ${\rm End}_k (\ddot N)$ is $0$-split. In this case, the embedded algebras $(k G)_\gamma$ and $(k G')_{\gamma'}$ associated to $P_\gamma$ and $P_{\gamma'}$ are source algebras (see \cite{T}) of the block algebras $k G b$ and $k G'b'$ respectively.

\medskip\noindent{\bf 14.}\quad
Clearly we have the equality
${\Bbb C}_0({\rm End}_k (\ddot N)\otimes_ k (k G')_{\gamma'})={\Bbb C}_0({\rm End}_k (\ddot N))\otimes_ k (k G')_{\gamma'}$ which induces a $k P$-interior algebra isomorphism  $${\Bbb H}_0({\rm End}_k (\ddot N)\otimes_ k(k G')_{\gamma'})\cong {\Bbb H}_0({\rm End}_k (\ddot N))\otimes_ k(k G')_{\gamma'}.\leqno 14.1$$
Since $P$ is normal in $G'$ and the $kP$-interior algebra ${\Bbb C}_0({\rm End}_k (\ddot N))$ is a local primitive $k P$-interior algebra, by Lemma 11 the point $\ddot\gamma$ only contains the
identity element of ${\rm End}_k (\ddot N)\otimes_ k(k G')_
{\gamma'}$. Then by Isomorphisms 13.1 and 14.1 we have a $k P$-interior algebra isomorphism
$$(k G)_\gamma\cong {\Bbb H}_0({\rm End}_k (\ddot N))\otimes_ k(k G')_{\gamma'}.\leqno 14.2$$
We consider the ${\frak D} P$-modules ${\rm End}_k (\ddot N)$ and
${\rm End}_k (\ddot N)\otimes_ k (k G')_{\gamma'}$
determined respectively by their $\frak D$-algebra structures and the $P$-conjugations on them.
Since $(k G')_{\gamma'}$ has a $P\times P$-stable basis containing its unity (see \cite[Proposition 38.7]{T}), the ${\frak D} P$-module
${\rm End}_k (\ddot N)$ is a direct summand of the ${\frak D}P$-module ${\rm End}_k (\ddot N)\otimes_ k (k G')_{\gamma'}$. Since
the ${\frak D}({P}\times {P})$-module ${\rm End}_k (\ddot N)\otimes_ k (k G')_
{\gamma'}$ is $0$-split, so is the ${\frak D} P$-module
${\rm End}_k (\ddot N)$.

\medskip\noindent{\bf 15.}\quad
Since $(k G)_\gamma$ is a primitive $k P$-interior algebra, the isomorphism 14.2 forces that the $k P$-interior algebra ${\Bbb H}_0({\rm End}_k (\ddot N))$ is also a primitive $k P$-interior algebra.
Then by Lemma 12, the homology of the $\frak D$-module $\ddot N$ vanishes at all degree but degree $q$ and
the map
${\Bbb C}_0({\rm End}_k(\ddot N))\rightarrow {\rm End}_k({\Bbb H}_q(\ddot N))$ sending a chain map $f: \ddot N\rightarrow \ddot N$ onto the induced $k$-module homomorphism $f_q: {\Bbb H}_q(\ddot N)\rightarrow {\Bbb H}_q(\ddot N)$ induces a $k$-algebra isomorphism$${\Bbb H}_0({\rm End}_k(\ddot N))\cong {\rm End}_k({\Bbb H}_q(\ddot N));\leqno 15.1$$ moreover this $k$-algebra isomorphism actually is a $k P$-interior algebra isomorphism. So we have a $k P$-interior algebra
isomorphism  $(k G)_\gamma\cong{\rm End}_k({\Bbb H}_q(\ddot N))\otimes_ k(k G')_{\gamma'}$.
By \cite[Theorem 7.2]{P3}, ${\Bbb H}_q(\ddot N)$ is an endo-permutation $k P$-module with vertex $P$.

\medskip\noindent{\bf 16.}\quad
Clearly $\{{\rm id}_{\ddot M}\}$ is a point of $G\times G'$ on the ${\frak D}(G\times G')$-interior algebra ${\rm End}_k(\ddot M)$, where ${\rm id}_{\ddot M}$ is the identity map on $\ddot M$.
Since the ${\frak D}P$-module $\ddot N$ is a source of the ${\frak D}(G\times G')$-module $\ddot M$, the ${\frak D}P$-module $\ddot N$ corresponds to a pointed group $P_{\ddot\delta}$ on ${\rm End}_k(\ddot M)$ so that $P_{\ddot\delta}$ is a defect pointed group of the pointed group $(G\times G')_{\{{\rm id}_{\ddot M}\}}$; moreover the ${\frak D}P$-interior algebra $ {\rm End}_k(\ddot N)$ is an embedded algebra associated with $P_{\ddot\delta}$. We denote by ${\rm Ind}^{G\times G'}_P({\rm End}_k(\ddot N))$ the injective induction of the ${\frak D}P$-interior algebra ${\rm End}_k(\ddot N)$ from $P$ to $G\times G'$ (see \cite[12.2]{P3}), which is a ${\frak D}(G\times G')$-interior algebra.
Then by \cite[Corollary 14.11]{P3}, there is a ${\frak D}(G\times G')$-interior algebra embedding $$\ddot{\rm H}:{\rm End}_k(\ddot M)\rightarrow {\rm Ind}^{G\times G'}_P({\rm End}_k(\ddot N)).$$
On the other hand, by
\cite[2.6.5]{P3}, it is easily checked that there is a ${\frak D}(G\times G')$-interior algebra isomorphism ${\rm Ind}^{G\times G'}_P({\rm End}_k(\ddot N))\cong {\rm End}_k({\rm Ind}^{G\times G'}_P(\ddot N))$. Therefore we get a ${\frak D}(G\times G')$-interior algebra embedding ${\rm End}_k(\ddot M)\rightarrow{\rm End}_k({\rm Ind}^{G\times G'}_P(\ddot N))$. Then by \cite[Example 13.4]{T}, it is easily checked that $\ddot M$ is a direct summand of the ${\frak D}(G\times G')$-module ${\rm Ind}^{G\times G'}_P(\ddot N)$.

\medskip\noindent{\bf 17.}\quad Since the homology of the $\frak D$-module $\ddot N$ vanishes at all degree but degree $q$, so do the homology of the ${\frak D}(G\times G')$-module ${\rm Ind}^{G\times G'}_P(\ddot N)$ and the homology of the ${\frak D}(G\times G')$-module $\ddot M$; moreover we have ${\Bbb H}_q({\rm Ind}^{G\times G'}_P(\ddot N))={\rm Ind}^{G\times G'}_P({\Bbb H}_q(\ddot N))$.  Then by \cite[Theorem 3.1]{HS},
the map
$${\Bbb H}_0({\rm End}_k({\rm Ind}^{G\times G'}_P(\ddot N)))\rightarrow {\rm End}_k({\rm Ind}^{G\times G'}_P({\Bbb H}_q(\ddot N)))\leqno 17.1$$ sending the image in ${\Bbb H}_0({\rm End}_k({\rm Ind}^{G\times G'}_P(\ddot N)))$ of a chain map $f: {\rm Ind}^{G\times G'}_P(\ddot N)\rightarrow {\rm Ind}^{G\times G'}_P(\ddot N)$ onto the induced $k$-module homomorphism $f_q: {\rm Ind}^{G\times G'}_P({\Bbb H}_q(\ddot N))\rightarrow {\rm Ind}^{G\times G'}_P({\Bbb H}_q(\ddot N))$ is a $k$-linear isomorphism, which actually is a $k(G\times G')$-interior algebra isomorphism. Since the ${\frak D}(G\times G')$-module $\ddot M$ is a direct summand of ${\rm Ind}^{G\times G'}_P(\ddot N)$,  Isomorphism 17.1 induces a $k(G\times G')$-interior algebra isomorphism
$\Phi: {\Bbb H}_0({\rm End}_k(\ddot M))\rightarrow {\rm End}_k({\Bbb H}_q(\ddot M))$.

\medskip\noindent{\bf 18.}\quad We consider the $k P$-interior algebra homomorphism ${\Bbb C}_0({\rm End}_k(\ddot N))\rightarrow {\rm End}_k({\Bbb H}_q(\ddot N))$ obtained by composing the surjective homomorphism ${\Bbb C}_0({\rm End}_k(\ddot N))\rightarrow {\Bbb H}_0({\rm End}_k(\ddot N))$ and Isomorphism 15.1.  Since ${\Bbb C}_0({\rm Ind}^{G\times G'}_P({\rm End}_k(\ddot N)))={\rm Ind}^{G\times G'}_P({\Bbb C}_0({\rm End}_k(\ddot N)))$,
the homomorphism ${\Bbb C}_0({\rm End}_k(\ddot N))\rightarrow {\rm End}_k({\Bbb H}_q(\ddot N))$ induces a surjecitve $k(G\times G')$-interior algebra homomorphism $${\Bbb C}_0({\rm Ind}^{G\times G'}_P({\rm End}_k(\ddot N)))\rightarrow {\rm Ind}^{G\times G'}_P({\rm End}_k({\Bbb H}_q(\ddot N)))$$ with the kernel ${\Bbb B}_0({\rm Ind}^{G\times G'}_P({\rm End}_k(\ddot N)))$, which induces a $k(G\times G')$-interior algebra isomorphism $$\Psi: {\Bbb H}_0({\rm Ind}^{G\times G'}_P({\rm End}_k(\ddot N)))\cong {\rm Ind}^{G\times G'}_P({\rm End}_k({\Bbb H}_q(\ddot N))).$$

\medskip\noindent{\bf 19.}\quad Clearly the embedding $\ddot{\rm H}$ induces a $k(G\times G')$-interior algebra embedding $${\Bbb C}_0({\rm End}_k(\ddot M))\rightarrow {\Bbb C}_0({\rm Ind}^{G\times G'}_P({\rm End}_k(\ddot N)))$$ which maps ${\Bbb B}_0({\rm End}_k(\ddot M))$
into ${\Bbb B}_0({\rm Ind}^{G\times G'}_P({\rm End}_k(\ddot N)))$. Thus the embedding $\ddot{\rm H}$ induces
a $k(G\times G')$-interior algebra embedding ${\Bbb H}_0(\ddot{\rm H}):{\Bbb H}_0({\rm End}_k(\ddot M))\rightarrow {\Bbb H}_0({\rm Ind}^{G\times G'}_P({\rm End}_k(\ddot N)))$.
We set $${\rm H}=\Psi\circ{\Bbb H}_0(\ddot{\rm H})\circ\Phi^{-1}.$$ We denote by $\phi$ the composition of the canonical surjective homomorphism ${\Bbb C}_0({\rm End}_k(\ddot M))\rightarrow {\Bbb H}_0({\rm End}_k(\ddot M))$ with the isomorphism $\Phi$ and by $\psi$ the compostion of the canonical surjective homomorphism ${\Bbb C}_0({\rm Ind}^{G\times G'}_P({\rm End}_k(\ddot N)))\rightarrow {\Bbb H}_0({\rm Ind}^{G\times G'}_P({\rm End}_k(\ddot N)))$ with the isomorphism $\Psi$.
Then we have the following commutative diagram
$$\matrix{{\Bbb C}_0({\rm End}_k(\ddot M))&\longrightarrow & {\Bbb C}_0({\rm Ind}^{G\times G'}_P({\rm End}_k(\ddot N)))\cr
\phi\Big{\downarrow}&&\Big{\downarrow}\psi\cr
{\rm End}_k({\Bbb H}_q(\ddot M))&\buildrel{\rm H}\over{\longrightarrow}  &{\rm Ind}^{G\times G'}_P({\rm End}_k({\Bbb H}_q(\ddot N))).}\leqno 19.1$$

\medskip\noindent{\bf 20.}\quad Now we claim that the homomorphism $\phi$ induces a surjective $k G$-interior algebra homomorphism ${\Bbb C}_0({\rm End}_{k(1\times G')}(\ddot M))\rightarrow {\rm End}_{k(1\times G')}({\Bbb H}_q(\ddot M))$ with the kernel ${\Bbb B}_0({\rm End}_{k(1\times G')}(\ddot M))$. Since we have the commutative diagram 19.1, it suffices to show that the homomorphism $\psi$ induces a surjective $k G$-interior algebra homomorphism ${\Bbb C}_0({\rm Ind}^{G\times G'}_P({\rm End}_k(\ddot N))^{1\times G'})\rightarrow {\rm Ind}^{G\times G'}_P({\rm End}_k({\Bbb H}_q(\ddot N)))^{1\times G'}$ with the kernel ${\Bbb B}_0({\rm Ind}^{G\times G'}_P({\rm End}_k(\ddot N))^{1\times G'})$. By \cite[Theorem 15.4]{P3}, there are a unique ${\frak D}G$-interior algebra isomorphism  $${\rm H}_{G,\,G'}^{\ddot N}:{\rm Ind}^{G\times G'}_P({\rm End}_k(\ddot N))^{1\times G'}\cong{\rm Ind}_P^G({\rm End}_k(\ddot N)\otimes_k k G') $$ mapping ${\rm Tr}^{1\times G'}_{1\times 1}(x\otimes a\otimes 1)$ onto $1\otimes (a\otimes x^{-1})\otimes 1$ for any $a\in {\rm End}_k(\ddot N)$ and any $x\in G'$. This isomorphism ${\rm H}_{G,\,G'}^{\ddot N}$ induces a $k G$-interior algebra isomorphism $${\Bbb C}_0({\rm H}_{G,\,G'}^{\ddot N}):{\rm Ind}^{G\times G'}_P({\Bbb C}_0({\rm End}_k(\ddot N)))^{1\times G'}\cong{\rm Ind}_P^G({\Bbb C}_0({\rm End}_k(\ddot N))\otimes_k k G').$$

\medskip\noindent{\bf 21.}\quad By \cite[Theorem 4.4]{P3}, there is a unqiue $kG$-interior algebra isomorphism $${\rm H}_{G,\,G'}^{{\Bbb H}_q(\ddot N)}:{\rm Ind}^{G\times G'}_P({\rm End}_k({\Bbb H}_q(\ddot N)))^{1\times G'}\cong{\rm Ind}_P^G({\rm End}_k({\Bbb H}_q(\ddot N))\otimes_k  k G')$$mapping ${\rm Tr}^{1\times G'}_{1\times 1}(x\otimes a\otimes 1)$ onto $1\otimes (a\otimes x^{-1})\otimes 1$ for any $a\in {\rm End}_k({\Bbb H}_q(\ddot N))$ and any $x\in G'$. Clearly the $k P$-interior algebra homomorphism ${\Bbb C}_0({\rm End}_k(\ddot N))\rightarrow {\rm End}_k({\Bbb H}_q(\ddot N))$ (see Paragraph 17) induces a surjective $k G$-interior algebra homomorphism $${\rm Ind}_P^G({\Bbb C}_0({\rm End}_k(\ddot N))\otimes_k k G')\rightarrow {\rm Ind}_P^G({\rm End}_k({\Bbb H}_q(\ddot N))\otimes_k k G')$$
with the kernel ${\rm Ind}_P^G({\Bbb B}_0({\rm End}_k(\ddot N))\otimes k G')={\Bbb B}_0({\rm Ind}_P^G({\rm End}_k(\ddot N)\otimes k G')$,  and this induced $k G$-interior algebra homomorphism makes the following diagram commutative
$$\matrix{{\rm Ind}^{G\times G'}_P({\Bbb C}_0({\rm End}_k(\ddot N)))^{1\times G'}&\buildrel {\Bbb C}_0({\rm H}_{G,\,G'}^{\ddot N})\over{\longrightarrow }& {\rm Ind}_P^G({\Bbb C}_0({\rm End}_k(\ddot N))\otimes_k k G')\cr
\psi\Big{\downarrow}&&\Big{\downarrow}\cr
{\rm Ind}^{G\times G'}_P({\rm End}_k({\Bbb H}_q(\ddot N)))^{1\times G'}&\buildrel {\rm H}_{G,\,G'}^{{\Bbb H}_q(\ddot N)}\over{\longrightarrow}  &{\rm Ind}_P^G({\rm End}_k({\Bbb H}_q(\ddot N))\otimes_k k G').}$$ Therefore the $k G$-interior algebra homomorphism $${\Bbb C}_0({\rm Ind}^{G\times G'}_P({\rm End}_k(\ddot N)))^{1\times G'}\rightarrow {\rm Ind}^{G\times G'}_P({\rm End}_k({\Bbb H}_q(\ddot N)))^{1\times G'}$$ induced by $\psi$ is surjective with the kernel ${\Bbb B}_0({\rm Ind}^{G\times G'}_P({\rm End}_k(\ddot N))^{1\times G'})$.  Then the claim is done.

\medskip\noindent{\bf 22.}\quad In particular, we have a $k G$-interior algebra isomorphism $${\Bbb H}_0({\rm End}_{k(1\times G')}(\ddot M))\cong {\rm End}_{k(1\times G')}({\Bbb H}_q(\ddot M)).$$ On the other hand, since the ${\frak D}(G\times G')$-module $\ddot M$ induces a Rickard equivalence between $k Gb$ and $kG'b'$, by \cite[Theorem 18.4]{P3} we have a $k G$-interior algebra isomorphism $k G b\cong {\Bbb H}_0({\rm End}_{k(1\times G')}(\ddot M)).$ Therefore we have a $k G$-interior algebra isomorphism $$k G b\cong {\rm End}_{k(1\times G')}({\Bbb H}_q(\ddot M)).$$ Then by \cite[Proposition 6.5]{P3} the $k(G\times G')$-module ${\Bbb H}_q(\ddot M)$ induces a Morita equivalence between $k Gb$ and $k G' b'$.
We claim that this Morita equivalence is basic in the sense of Puig in \cite{P3}.
There is a $k(G\times G')$-interior algebra isomorphism $${\rm Ind}^{G\times G'}_P({\rm End}_k({\Bbb H}_q(\ddot N)))\cong {\rm End}_k({\rm Ind}^{G\times G'}_P({\Bbb H}_q(\ddot N))).$$By composing $\rm H$ with this $k(G\times G')$-interior algebra isomorphism, we get a $k(G\times G')$-interior algebra embedding $ {\rm End}_{k}({\Bbb H}_q(\ddot M))\rightarrow {\rm End}_k({\rm Ind}^{G\times G'}_P({\Bbb H}_q(\ddot N)))$. Then by \cite[Example 13.4]{P3},
$\ddot M$ is a direct summand of the $k(G\times G')$-module ${\rm Ind}^{G\times G'}_P({\Bbb H}_q(\ddot N))$. Since the ${\frak D}P$-module ${\Bbb H}_q(\ddot N)$ is a direct summand of the restriction ${\rm Res}^{G\times G'}_{\ddot P}(M)$ and has vertex $\Delta(P)$, the ${\frak D}P$-module ${\Bbb H}_q(\ddot N)$ is  a source of the ${\frak D}(G\times G')$-module $M$. The claim is done.

\bigskip Now it remains to prove Statement 8.1 from Statement 8.2.  We continue to keep the notation in Theorem 8, and assume that Statement 8.2 holds.

\medskip\noindent{\bf 23.}\quad Since the ${\frak D}\ddot P$-module ${\rm End}_k(\ddot N)$ is $0$-split, the ${\frak D}\ddot P$-module ${\rm End}_k(\ddot N)$ is a direct sum of a contractible ${\frak D}\ddot P$-module and the $k \ddot P$-module ${\Bbb H}_0({\rm End}_k(\ddot N))$ determined by the $\ddot P$-conjugation. That is to say, there is a contractible ${\frak D}\ddot P$-module $C$ such that we have a direct sum decomposition ${\rm End}_k(\ddot N)={\Bbb H}_0({\rm End}_k(\ddot N))\oplus C$, from which, we easily conclude a direct sum decomposition $$({\rm End}_k(\ddot N))^{\ddot P}_R={\Bbb H}_0({\rm End}_k(\ddot N))^{\ddot P}_R\oplus C^{\ddot P}_R,\leqno 23.1$$  for any subgroup $R$ of $\ddot P$.  Since the ${\frak D}\ddot P$-module $C$ is contractible, so is the $\frak D$-module $ C^{\ddot P}_R$. Thus the $\frak D$-module $({\rm End}_k(\ddot N))^{\ddot P}_R$ is $0$-split and we have ${\Bbb C}_0(C^{\ddot P}_R)={\Bbb B}_0(C^{\ddot P}_R)={\Bbb B}_0(({\rm End}_k(\ddot N))^{\ddot P}_R)$. From the decomposition 23.1,  we get a new direct sum decomposition
\begin{center} ${\Bbb C}_0(({\rm End}_k(\ddot N))^{\ddot P}_R)={\Bbb H}_0({\rm End}_k(\ddot N))^{\ddot P}_R\oplus {\Bbb C}_0(C^{\ddot P}_R),$  \end{center}for any subgroup $R$ of $\ddot P$. This new decomposition implies that
the inclusion map ${\rm End}_{k \ddot P}(\ddot N)\subset {\rm End}_{k}(\ddot N)$ induces a surjective $k$-algebra homomorphism
$${\Bbb C}_0({\rm End}_{k \ddot P}(\ddot N))\rightarrow {\Bbb H}_0({\rm End}_{k}(\ddot N))^{\ddot P}\leqno 23.2$$ with the kernel ${\Bbb B}_0({\rm End}_{k{\ddot P}}(\ddot N))$, which induces a surjective $k$-algebra homomorphism
$${\Bbb C}_0({\rm End}_{k \ddot P}(\ddot N))\rightarrow {\Bbb H}_0({\rm End}_{k}(\ddot N))({\ddot P})$$ with the kernel ${\Bbb B}_0({\rm End}_{k{\ddot P}}(\ddot N))+\sum_R({\Bbb C}_0({\rm End}_k(\ddot N)))^{\ddot P}_R,$ where $R$ runs over proper subgroups of $\ddot P$.

\medskip\noindent{\bf 24.}\quad Since $\ddot M$ is a noncontractible ${\frak D}(G\times G')$-module, so is the ${\frak D}\ddot P$-module $\ddot N$. So ${\Bbb C}_0({\rm End}_{k\ddot P}(\ddot N))$ is a pimitive $k$-algebra with the identity map outside ${\Bbb B}_0({\rm End}_{k\ddot P}(\ddot N))$.
Then it follows from Homomorphism 23.2 that ${\Bbb H}_0({\rm End}_{k}(\ddot N))$ is a primitive $k \ddot P$-interior algebra. Then by Lemma 12, there is a $k \ddot P$-interior algebra isomorphism ${\Bbb H}_0({\rm End}_k(\ddot N))\cong {\rm End}_k({\Bbb H}_q(\ddot N));$ moreover, the $k\ddot P$-module ${\Bbb H}_q(\ddot N)$ is indecomposable.  We claim that  the $k\ddot P$-module ${\Bbb H}_q(\ddot N)$ has vertex $\ddot P$. Otherwise, by \cite[Proposition 18.11]{T} we have ${\rm End}_k({\Bbb H}_q(\ddot N))(\ddot P)=0$ and then ${\Bbb H}_0({\rm End}_{k}(\ddot N))({\ddot P})=0$; thus by Rosenberg's Lemma, ${\rm id}_{\ddot N}$ belongs to either ${\Bbb B}_0({\rm End}_{k{\ddot P}}(\ddot N))$ or $({\Bbb C}_0({\rm End}_k(\ddot N)))^{\ddot P}_R$ for some proper subgroup $R$ of $\ddot P$; this contradicts with the ${\frak D}\ddot P$-module $\ddot N$ being a source of the noncontractible ${\frak D}(G\times G')$-module $\ddot M$. The claim is done.

\medskip\noindent{\bf 25.}\quad Next we claim that the $k \ddot P$-module ${\Bbb H}_q(\ddot N)$ is a source of the $k(G\times G')$-module ${\Bbb H}_q(\ddot M)$.
Since the ${\frak D}\ddot P$-module $\ddot N$ is a direct summand of the restricton ${\rm Res}^{G\times G'}_{\ddot P}(\ddot M)$,  the $kP$-module
${\Bbb H}_q(\ddot N)$ is a direct summand of the restriction ${\rm Res}^{G\times G'}_{\ddot P}({\Bbb H}_q(\ddot M))$. Since the ${\frak D}\ddot P$-module $\ddot N$ is a source of the ${\frak D}(G\times G')$-module $\ddot M$, $\ddot M$ is a direct summand of the induced ${\frak D}(G\times G')$-module ${\rm Ind}^{G\times G'}_{\ddot P}(\ddot N)$ (see Paragraph 16), thus ${\Bbb H}_q(\ddot M)$ is a direct summand of the $k(G\times G')$-module ${\Bbb H}_q({\rm Ind}^{G\times G'}_{\ddot P}(\ddot N))$. Clearly there is an obvious $k(G\times G')$-module isomorphism ${\Bbb H}_q({\rm Ind}^{G\times G'}_{\ddot P}(\ddot N))\cong {\rm Ind}^{G\times G'}_{\ddot P}({\Bbb H}_q(\ddot N))$. So ${\Bbb H}_q(\ddot M)$ is a direct summand of the $k(G\times G')$-module ${\rm Ind}^{G\times G'}_{\ddot P}({\Bbb H}_q(\ddot N))$. Since the indecomposable $k\ddot P$-module ${\Bbb H}_q(\ddot N)$ has vertex $\ddot P$, the $k\ddot P$-module ${\Bbb H}_q(\ddot N)$ is a source of the $k(G\times G')$-module ${\Bbb H}_q(\ddot M)$.

\medskip
Finally, since the module ${\Bbb H}_q(\ddot M)$ induces a basic Morita equivalence between $k G b$ and $k G'b'$, by \cite[Corollary 7.4]{P3} the groups $\ddot P$ and $P$ have the same order.

\end{document}